\documentclass[a4paper,12pt]{article}
\usepackage{makeidx}
\usepackage[T1]{fontenc}
\usepackage{amsmath,amscd,amsthm}
\usepackage{amssymb}
\usepackage{tabularx}
\usepackage{amssymb,eucal,bezier,graphicx,}
\usepackage{times,amssymb}
\usepackage[ansinew]{inputenc}
\usepackage{multicol}

\newtheorem{thm}{Theorem}[section]

\newtheorem{lemma}[thm]{Lemma}

\newtheorem{defi}[thm]{Definition}

\begin{document}

%AAAAAAAAAAAAAAAAAAAAAAAAAAAAAAAAAAAAAAAAAAAAAAAA

\renewcommand{\figurename}{Figure} %Cambia la palabra "Figura" por "Figure"

\renewcommand\thefigure{\arabic{section}.\arabic{figure}} % Genera numeración X.Y = X es la sección Y número de figura independientemente de cada sección

\numberwithin{figure}{section} %Hace que la primera figura de cada sección X sea X.1 = X es la sección

%\renewcommand{\thebibliography}{References}

%\addto\captionsspanish{%
%\def\bibname{References}%
%}

%AAAAAAAAAAAAAAAAAAAAAAAAAAAAAAAAAAAAAAAAAAAAAAAAAAA

%\begin{center} {\large \bf Associating certain graphs with graphicable algebras}
%\end{center}

\begin{center} {\large \bf Certain particular families of graphicable algebras}
\end{center}

\begin{center}
{\bf Juan N\'u\~nez, María Luisa Rodríguez-Arévalo and María Trinidad Villar}
\end{center}

\begin{center}
{\small Dpto. Geometr\'{\i}a y Topolog\'{\i}a. \\
\small Facultad de Matem\'aticas. Universidad de Sevilla. \\
\small Apdo. 1160. 41080-Sevilla, Spain.\\
\small jnvaldes@us.es \,\, ml.rodriguezarevalo@gmail.com \,\, villar@us.es }
\end{center}

\vspace{0.4cm}

\begin{center}
{\bf Abstract}

\end{center}

\begin{quotation}
\noindent In this paper, we introduce some particular families of graphicable algebras obtained by following a relatively new line of research, initiated previously by some of the authors. It consists of the use of certain objects of Discrete Mathematics, mainly graphs and digraphs, to facilitate the study of graphicable algebras, which are a subset of evolution algebras.
\end{quotation}

\vspace{0.1cm}


\vspace{0.1cm}


\section*{Introduction}

The main goal of this paper is to advance in the research of a novel mathematical topic emerged not long ago, the {\em evolution algebras} in general, and the {\em graphicable algebras} (a subset of them) in particular, in order to obtain new results starting from those by Tian (see \cite{book, Petr}) and others already obtained by some of us in a previous paper \cite{AMC}, also relying in conceptually separated tools from them, such as graphs and digraphs. Concretely, our goal is to find some particular types of graphicable algebras associated with well-known types of graphs.

The motivation to deal with evolution algebras in general and graphicable algebras in particular is due to the fact that at present, the study of these algebras is very booming, due to the numerous connections between them and many other branches of Mathematics, such as Graph Theory, Group Theory, Markov processes, dynamic systems and the Theory of Knots, among others. Furthermore, they are also related to other sciences, including Biology, since Mendelian genetics is precisely which originated them.

For this study we have used in this work certain objects of Discrete Mathematics, particularly graphs and digraphs, thus linking two branches of the Mathematic between which at first appears to be no connection: Graph Theory and algebras of evolution theory, so that the attainment of new properties of each of them allows to achieve certain advances in the study of the other.

The relationship between both subjects is based on the representation of each evolution algebra by a certain graph, the study of the properties of these graphs and its later translation to the field of evolution algebras. Finally, in a similar way, the reverse process is also attempted, that is, to use the algebraic properties of evolution algebras to obtain properties for some types of graphs and advance, therefore, in the study of these objects. This research, which could be considered novel,  was somehow used by Tian himself in Chapter 6 of his book \cite {book} as a possibility to progress in the knowledge of these algebras.

Concretely, in this paper, we start from well-known types of graphs and introduce their corresponding associated graphicable algebras. Among these so-defined algebras we deal in Section 2 with complete tripartite and n-partite, star, friendship, wheel and snark graphicable algebras and set certain links of inclusion among them (see Theorem \ref{Main}). As it had already been proved in \cite{AMC}, a subalgebra of a graphicable algebra is not always graphicable. Our study shows the existence of families of graphicable algebras in which this property is hereditary. Finally, other particular graphicable algebras associated with generalized Petersen graphs are also studied.

\section{Preliminaries}

In this section we indicate some basic concepts related to both Evolution Algebras and Discrete Mathematics, more specifically to Graph Theory, to be used in this paper. For a more general overview of both theories, \cite {book, Petr} for the first and \cite {clark} for the second, among others, are available. In this paper, all the evolution algebras will be considered finite.

{\em Evolution algebras}, which were firstly introduced  by J. P. Tian, and then jointly presented with Vojtechovsky in 2006 \cite{Petr}, and later appeared as a book by Tian in 2008 \cite{book}, are those in which the relationships between their generators $\{e_1, \ldots, e_n\}$ are given by
$$
\left\{
\begin{array}{l}
e_i  \cdot e_j = 0, \quad \forall i \neq j, \\
e_i  \cdot e_i = \displaystyle \sum_{j = 1}^v a_{ji} \, e_j.
\end{array}
\right.
$$

%\vspace{0.15cm}

         Tian and Vojtechovsky in \cite{Petr} proved that these algebras are commutative and flexible, in general, and non associative.

%\vspace{0.15cm}

In this paper we will deal with a subset of these algebras, namely graphicable algebras, introduced by Tian in \cite{book} in 2008.

%\vspace{0.15cm}

An {\em $n$-dimensional graphicable algebra} is a commutative, non associative algebra, with a set of generators $V = \{e_1, e_2, \dots, e_n \} $ endowed with relations
$$\left\{
\begin{array}{l}
e_i^2 = \displaystyle \sum_{k = 1}^n e_k, \quad 1 \leq i \leq k, \\
e_i \cdot e_j = 0, \quad i \neq j, \quad 1 \leq i, j \leq k. \\
\end{array}
\right.
\newline
$$

%\vspace{0.15cm}

Thus, it is obvious that a graphicable algebra is an evolution algebra, although the converse is not true in general.
%\vspace{0.15cm}

Moving on now Graph Theory, the most basic concepts, as simple graphs, pseudographs, directed graphs, edges, vertices, adjacency and incidence, adjacency matrix
and paths between two vertices of a graph are assumed to be known (\cite {clark} can be consulted for details).

A {\em morphism} between two simple graphs $G_1 = (V_1, E_1)$ and $G_2 = (V_2, E_2)$ is a map $ \phi $: $ V_ {1}$ $\rightarrow $ $V_ {2}$ such that if $u$ and $v$ are adjacent vertices in $G_ {1}$, then $ \phi (u)$ and $ \phi (v)$ are adjacent in $G_ {2}.$ Similarly, one can also define the concept of morphism between directed pseudo-graphs: a map verifying that if $(u, v)$ $\epsilon$ $E_{1}$, then $(\phi(u), \phi(v))$ $\epsilon$ $E_{2}$, taking into account the order of the pair of vertices defining the edge.

With respect to the association between graphicable algebras and graphs (where {\it graph} might be of any type above mentioned), Tian, in \cite{book}, showed how to associate a graph with a graphicable algebra; he gave the following

%\vspace{0.15cm}

\begin{defi}
Let $G = (V, E)$ be a graph, $V$ be the set of vertices and $E$ be the set of edges; It is defined the associated algebra with $G$ taking $V= \{e_1, e_2, \dots, e_r\}$ as the set of generators and $R$ as the set of relations of the algebra

$$
R = \left\{
\begin{array}{l}
e_i^2 = \displaystyle\sum_{e_k\in \Gamma(e_i)} e_k, \\
e_i \cdot e_j =0, \quad i \neq j,
\end{array}
\right.
\newline
$$
where $\Gamma(e_i)$ denotes the set of vertices adjacent to $e_i$.
\end{defi}

%\vspace{0.15cm}

Conversely, Tian also showed how to associate a graphicable algebra with a directed graph: he took the set of generators of algebra as the set of vertices and
as the set of edges those connecting the vertex $e_i$ with the vertices corresponding to generators appearing in the expression of $e_i^2$, for each generator $e_i$.

Tian himself, although by using a different notation to that followed here (he used $x_i$ for generators, instead of the notation that we use, $e_i$) introduced the graphicable algebras associated with cyclic graphs, paths, and complete graphs, referred to as cyclic algebras, path algebras and complete graphicable algebras, respectively.
%\vspace{0.15cm}

Later, N\'u\~nez, Silvero and Villar, following the study by Tian, obtained in \cite{AMC} graphicable algebras associated with cubic graphs $Q_6$ and $Q_8$ and with the Petersen and Heawood graphs. They also associated a graphicable algebra with complete bipartite graphs and wheel graphs. In that paper, they introduced the subfamily of $S$-graphicable algebras as those graphicable algebras verifying
$a_{ij}=a_{ji}$ for $i\neq j; \,\, i, j =1,\, 2, \dots, r,$ and $a_{ii}=0$ for every $i=1,\, 2,
\dots, r,$ where  $a_{ij}$ are the elements of the matrix associated to the algebra. Equivalently, $S$-graphicable algebras can be defined as those which have a simple graph associated  (see \cite{AMC}). In this paper we now introduce certain families of $S$-graphicable algebras starting from particular families of simple graphs.
Since all the graphicable algebras which appear from now here are $S$-graphicable, we will simply call them graphicable algebras.

%, for which the set of generators of the graphicable algebra associated to $W_n$ would be $ \{e_1, e_2, \dots, e_n, e_ {n +1} \}$, where the last generator is identified with the center of the algebra (by using this way of denoting the center of the wheel it can be seen that if $ e_{n +1} $ is suppressed in the law of algebra, the graphicable algebra associated with a cycle of $n$ vertices would be obtained).
%

\section{Particular types of graphicable algebras}

In this section we show our main results obtained following a procedure similar to that proposed by the last mentioned authors, but with the difficulties and problems emerged when considering particular graphs and specific families of graphs increasingly complicated in structure for having a higher number of edges and vertices. In particular, we show the $S$-graphicable algebras associated with $n$-partite graphs ($n\geq 3$), stars, wheels, friendship graphs, snark graphs and generalized Petersen graphs, respectively. Other specific families of graphs will be discussed in future work.

\subsection{Complete n-partite graphicable algebras ($n\geq 3$)}

It is already known (see \cite{AMC}) that given a complete bipartite graph, $K_{m,n}$, its associated graphicable algebra $A(K_{m,n})$ is defined as that one of dimension $m+n$ whose set of generators $\{e_1, e_2, \dots, e_{m+n}\}$ verifies the law given by the brackets

$$
\begin{array}{ll}
e_i^2 = \displaystyle{\sum_{k=1}^n e_{m +k}}, & \mbox{for}, 1 \leq i \leq m, \\
e_{m+i}^2  = \displaystyle{\sum_{k=1}^m e_{k}}, & \mbox{for}, 1 \leq i \leq n. \\
\end{array}
$$

So, similarly, we can define the graphicable algebra associated with a complete $n$-partite graph ($n\geq 3$), $K_{a_1,a_2,\ldots,a_n}$, that is, a graph where the $n$ subsets of  vertices have $a_1,a_2,\ldots,a_n$ elements, respectively.

%\begin{defi}
%It is defined the {\em tripartite graphicable algebra} $A(K_{m,n,r})$ as the graphicable algebra associated to a tripartite graph $K_{m,n,r}$ of $m+n+r$ vertices, that is, as the  graphicable algebra of dimension $m+n+r$, whose law with respect to the set of generators $
%\{e_1,\ldots,e_m,e_{m+1},e_{m+2},\ldots,e_{m+n},e_{m+n+1},\ldots,e_{m+n+r},\}$ is given by
%$$
%\begin{array}{ll}
%  e_i^2 = \displaystyle\sum_{k=1}^{n+r} e_{k+m}, & \mbox{ if } 1\leq i \leq m \\
%  e_i^2 = \displaystyle\sum_{k=1}^{m} e_{k} + \sum_{k=1}^{r} e_{k+m+n} & \mbox{ if } m < i \leq m+n \\
%  e_i^2 = \displaystyle\sum_{k=1}^{m+n} e_{k} & \mbox{ if } m+n < i \leq m+n+r.
%\end{array}
%$$
%\end{defi}

%Next, by using this procedure, we can generalize the previous definition for the case $n$-partite.

\begin{defi}
It is defined the {\em $n$-partite graphicable algebra} $A(K_{a_1,a_2,\ldots,a_n})$ ($n\geq 2$) as the graphicable algebra associated to a $n$-partite graph with $a_1 + a_2 + \ldots + a_n$ vertices, $K_{a_1,a_2,\ldots,a_n}$. This algebra is of dimension $a_1+\ldots+a_n$, with a set of generators
$\{e_1,\ldots,e_{a_1},\ldots,e_{a_1+a_2},\ldots,e_{a_1+\ldots+a_n},\}$ and law

\vspace{0.15cm}

$$
\begin{array}{ll}
  e_i^2 = \displaystyle\sum_{k=1}^{a_2+\ldots+a_n} e_{k+a_1}, & \mbox{ if } 1\leq i \leq a_1, \\
  e_i^2 = \displaystyle\sum_{k=1}^{a_1} e_{k} + \sum_{k=1}^{a_3+\ldots+a_n} e_{k+a_1+a_2} & \mbox{ if } a_1 < i \leq a_1+a_2, \\
  e_i^2 = \displaystyle\sum_{k=1}^{a_1+a_2} e_{k} + \sum_{k=1}^{a_3+\ldots+a_n}e_{k+a_1+a_2+a_3} & \mbox{ if } a_1+a_2 < i \leq a_1+a_2+a_3, \\
  \qquad \qquad  \qquad \quad \vdots & \qquad \qquad \qquad \quad \vdots \\
  e_i^2 = \displaystyle\sum_{k=1}^{a_1+\ldots+a_{n-1}} e_{k} & \mbox{ if } a_1+\ldots+a_{n-1} < i \leq a_1+\ldots+a_n.
\end{array}
$$
\end{defi}

\vspace{0.15cm}

\subsection{Star graphicable algebras}

A {\em star graph}, or simply a {\em star}, $S_n$, is the complete bipartite graph $K_{1,n}$, that is, a tree with one internal node and $n\geq 2$ leaves (but, no internal nodes and $n + 1$ leaves when $n = 1$). Other authors define $S_n$ to be the tree of order $n$ with maximum diameter $2$. From now on, the first criterion in the definition
of $S_n$ is adopted. Note that stars may also be described as the only connected graphs in which at most one vertex has degree greater than one. As a particular case, a star with 3 edges is called a {\em claw}.

\vspace{0.15cm}

\begin{figure}[h!]
\centering
\includegraphics[width=6cm]{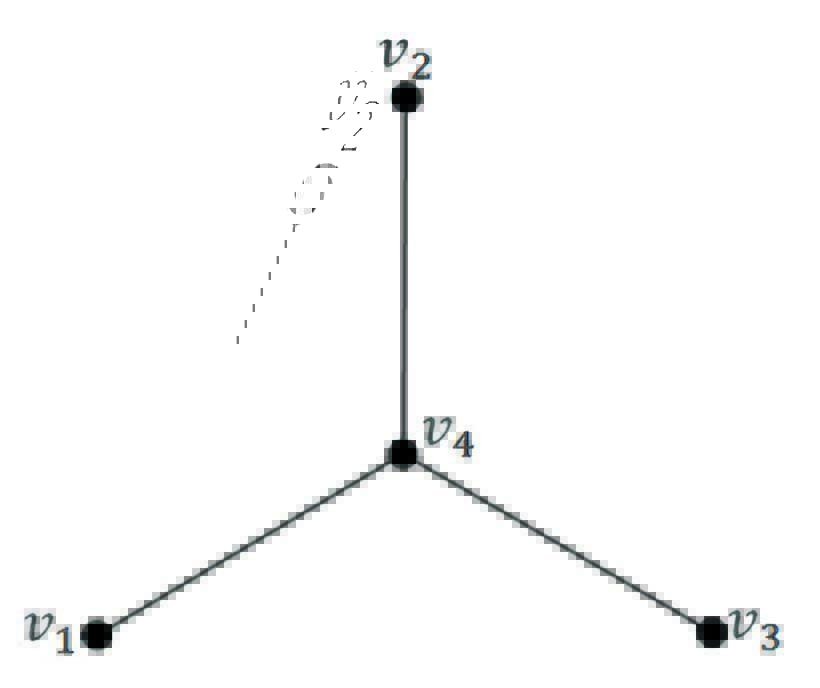}
\caption{Star graph $S_3$.}
\label{Figura1}
\end{figure}

Our objective is to associate this family of graphs with a family of graphicable algebras. To do this, as the simplest case $n=2$ is trivial, we will begin with the star which has three leaves and then try to generalize to the case of $n$ leaves by induction. Therefore, we consider the following star graph $S_3$ of Figure \ref{Figura1}.

\newpage

According to our procedure, the graphicable algebra associated to this graph has the generators $\{e_1,e_2,e_3,e_4\}$ and satisfies the law
$$
e_1^2= e_4, \quad e_2^2=e_4, \quad e_3^2=e_4, \quad e_4^2=e_1 + e_2 + e_3,
$$

%Similarly, the graphicable algebra associated to the star $S_4$ (see Figure \ref{Figura2})
%
%\begin{figure}[h!]
%\centering
%\includegraphics[width=4.9cm]{S4}
%\caption{Star graph $S_4$.}
%\label{Figura2}
%\end{figure}
%
%has the generators $\{e_1,e_2,e_3,e_4,e_5\}$ and its law is given by
%$$
% e_1^2= e_5, \quad  e_2^2=e_5, \quad   e_3^2=e_5, \quad   e_4^2=e_5, \quad e_5^2=e_1 + e_2 + e_3 + e_4.
%$$

%\begin{figure}[h!]
%\centering
%\includegraphics[width=4.9cm]{S4}
%\caption{Grafo estrella $S_4$.}
%\end{figure}

\noindent and the same could be done with stars $S_4$ and so on.
So, as a consequence, we can give the following

\begin{defi}
The graphicable algebra $A(S_n$) associated to the star $S_n$ is called {\em star graphicable algebra}. With respect to the set of generators $\{e_1,\ldots,e_{n+1}\}$ (where $e_{n+1}$ is identified with the center vertex), the law of this algebra is given by
$$
\begin{array}{l}
  e_i^2=e_{n+1}, \quad \forall i\neq n+1, \\
  e_{n+1}^2=\displaystyle\sum_{i=1}^{n} e_i.
\end{array}
$$
\end{defi}

\subsection{Friendship graphicable algebras}

\vspace{0.30cm}

%\begin{enumerate}
%\item \textbf{2.- graphs de la amistad}
%\end{enumerate}

%\vspace{0.15cm}

%\begin{itemize}
%\item \textbf{Friendship graphs}
%\end{itemize}

The {\em friendship graph} $F_n$, also named {\em Dutch windmill graph or n-fan graph} is a planar simple graph with  $2n+1$ vertices and $3n$ edges.

\vspace{0.15cm}

$F_n$ can be constructed by joining $n$ copies of the cycle graph $C_3$ with a common vertex. Due to that construction, $F_n$ is isomorphic to the windmill graph $Wd(3,n)$. In particular, $F_2$ is isomorphic to the butterfly graph.

\vspace{0.15cm}

Let us recall that the {\em friendship Theorem}, by Paul Erdös, Renyi Alfréd and Vera T. Sós, in 1966 (see \cite{Erdos}), states that finite graphs with the property that every two vertices have exactly one common adjacent vertex are the friendship graphs. Informally, if a group of people has the property that each pair of them have exactly one friend in common, then there must be a person who is friend with everyone else.

\vspace{0.15cm}

\begin{figure}[h!]
\centering
\includegraphics[width=5.8cm]{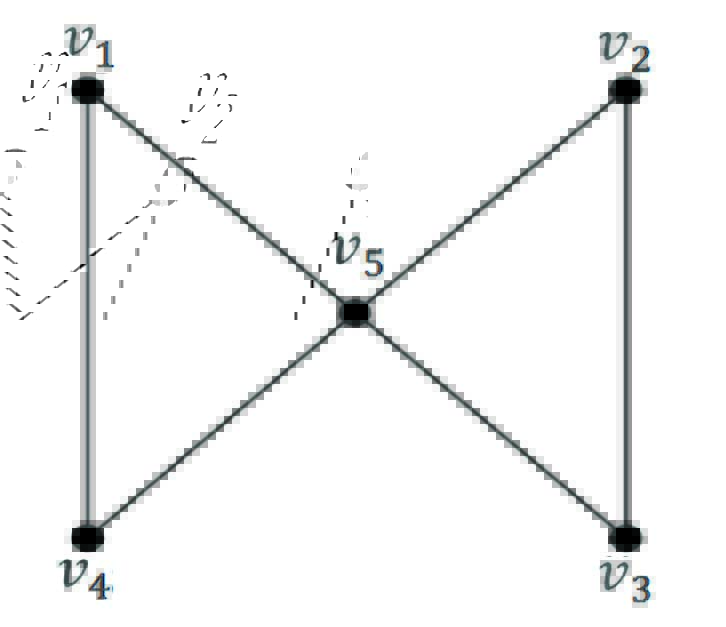}
\caption{Friendship graph $F_2$ or Butterfly graph.}
\label{Figura3}
\end{figure}

Let us now consider the friendship graph $F_2$ of Figure \ref{Figura3}

The associated graphicable algebra is the one of dimension 5 with  generators $\{e_1,e_2,e_3,e_4,e_5\}$ and the following law
$$
\begin{array}{l}
  e_1^2=e_4+e_5, \qquad e_2^2=e_3+e_5, \\
  e_3^2=e_2+e_5, \qquad e_4^2=e_1+e_5, \\
  e_5^2=e_1+e_2 + e_3 + e_4.
\end{array}
$$

\vspace{0.15cm}

%\newpage

Note that these relations can be written in the following way

\begin{itemize}
    \item ${\small {\it For \,\, i=1,\ldots,4, \,\, (mod.\,\,4),}}$
  \begin{itemize}
    \item ${\small {e_i^2=e_{i+1}+e_{5}}, \,\, \,\, \it if \,\, i \,\, is \,\, even},$
    \item ${\small {e_i^2=e_{i-1}+e_{5}}, \,\, \,\, \it if \,\, i \,\, is \,\, odd}.$
  \end{itemize}
  \item ${\small {\it For \,\, i=5, \,\, e_5^2=\displaystyle\sum_{i=1}^{4}e_i.}}$
\end{itemize}

%\vspace{0.15cm}
%
%For the case of graph $F_3$ we can obtain with a similar reasoning the corresponding associated graphicable algebras (see Figure \ref{Figura4})
%
%
%\begin{figure}[h!]
%\centering
%\includegraphics[width=6.5cm]{F3}
%\caption{friendship graph $F_3$.}
%\label{Figura4}
%\end{figure}
%
%
%\begin{itemize}
%    \item ${\small {\it For \,\, i=1,\ldots,6, \,\, (mod. \,\, 6),}}$
%  \begin{itemize}
%    \item ${\small {\it if \,\, i \,\, is \,\, even, \,\, e_i^2=e_{i+1}+e_{7}}},$
%    \item ${\small {\it if \,\, i \,\, is \,\, odd, \,\, e_i^2=e_{i-1}+e_{7}.}}$
%  \end{itemize}
%  \item ${\small {\it For \,\, i=7, \,\, e_7^2=\displaystyle\sum_{i=1}^{6}e_i.}}$
%\end{itemize}
%
%
%
%
%\begin{figure}[h!]\label{Figura5}
%\centering
%\includegraphics[width=3.5cm]{F4}
%\caption{friendship graph $F_4$.}
%\end{figure}
%
%
%\begin{itemize}
%  \item ${\small {\it For \,\, i=1, \,\, e_1^2=\displaystyle\sum_{i=2}^{9}e_i.}}$
%  \item ${\small {\it For \,\ i=2,\ldots,9,}}$
%  \begin{itemize}
%    \item ${\small {\it if \,\, i \,\, is \,\, even, \,\, e_i^2=e_1+e_{i+1} \,\, (mod. \, 5)}},$
%    \item ${\small {\it if \,\, i \,\, is \,\, odd, \,\, e_i^2=e_1+e_{i-1} \,\, (mod. \, 5).}}$
%  \end{itemize}
%\end{itemize}

\vspace{0.15cm}

Then, a similar procedure allows us to generalize these results obtaining the corresponding graphicable algebras associated with the family $F_n$ of friendship graphs.
Thus, we can give the following

\begin{defi}

The {\em friendship graphicable algebras} $A(F_n)$ are those associated to friendship graphs $F_n$ (with center vertex $v_{2n+1}$), that is, those of dimension $2n+1$, with generators $\{e_1,\ldots,e_{2n+1}\}$ and law given by

\vspace{0.15cm}

\begin{itemize}
    \item ${\small {\it For \,\, i=1,\ldots,2n, \,\, (mod. \,\, 2n),}}$
  \begin{itemize}
    \item ${\small {e_i^2=e_{i+1}+e_{2n+1}}, \quad \it if \,\, i \,\, is \,\, even},$
    \item ${\small {e_i^2=e_{i-1}+e_{2n+1}}, \quad \it if \,\, i \,\, is \,\, odd}.$
  \end{itemize}
  \item ${\small {\it For \,\, i=2n+1, \,\, e_{2n+1}^2=\displaystyle\sum_{i=1}^{2n}e_i.}}$
\end{itemize}
\end{defi}

\subsection{Wheel graphicable algebras}

Let us recall that a wheel graph $(W_n)$ is a graph with $n$ vertices ($n \geq 4$), formed by connecting a single vertex to all vertices of an $(n-1)$-cycle.

By using the same reasoning as before, we give the following

\begin{defi}
A {\em wheel graphicable algebra} $A(W_n)$ is that associated to a wheel graph $W_n$, that is, the one defined by the law

\begin{itemize}
  \item {$\small {\it For \,\, i=1,\ldots,n-1, \quad \displaystyle e_i^2=e_{i-1}+e_{i+1}+e_n, \,\, (mod. \,\, n-1).}$}
  \item ${\small {\it For \,\, i=n, \quad e_n^2=\displaystyle\sum_{i=1}^{n-1}e_i,}}$
\end{itemize}
\noindent with respect to the set of generators $\{e_1,\ldots,e_n\}$, where $e_n$ is the center of the graph $W_n$.
\end{defi}

\vspace{0.15cm}

Now, in the next subsection we are going to compare star, friendship and wheel graphicable algebras depending on the corresponding star, friendship and wheel graphs.

\vspace{0.15cm}

\subsection{Relationships among star, friendship and wheel graphicable algebras}

%Let consider the following graphs
%
%\begin{figure}[h!]
%\centering
%\includegraphics[width=13cm]{W7F3}
%\caption{Grafo rueda $W_7$ y grafo de la amistad $F_3$.}
%\end{figure}

We think that it is convenient to remark that the interest of the study of the latter families of graphs considered is the fact that some of its properties can be easily translated into the language of graphicable algebras, as discussed below.

\vspace{0.15cm}

In the first place, it is easy to realize that there is a relationship between friendship, star and wheel graphs. Concretely, each star graph is a subgraph of a given friendship graph. Similarly, each friendship is in turn a subgraph of a certain wheel. Explicitly, $S_n \subset F_{[n/2]} \subset W_n$. This allows us to deduce that a  friendship graphicable algebra can be constructed from the corresponding wheel graphicable algebra and vice versa. The same goes for star and the two other families.
Indeed, a friendship graph can be obtained from a wheel by suppressing edges and vertices in an appropriate way.

%between the vertex with subindex $i$ (according to the fixed notation) and the vertex with subindex $i + 1$, where $i$ takes odd values greater than $1$ in this notation.

\vspace{0.15cm}

Translating this fact into the language of graphicable algebras, it allows us to set the following

\begin{lemma}
In the expressions of $e_i^2$, with $1\leq i \leq 2n$, in the definition of friendship graphicable algebras $A(F_n)$ will appear a term less than in the corresponding wheel graphicable algebra $A(W_{2n+1})$. This term is $e_{i+1}$, if $i$ is odd or $e_{i-1}$ if $i$ is even (mod. $2n$). \hfill  $\Box$
\end{lemma}

\vspace{0.15cm}

In the same way, star graphs are also related with friendships and wheel graphs. In the first place, it is also easy to see that a friendship graph is obtained from a star graph of odd order by adding edges between the vertices of indexes $i$ and $i + 1$, where $i$ takes even values in the adopted notation.

\vspace{0.15cm}

Then, in the language of graphicable algebras, the following result is satisfied

\begin{lemma}
The expressions of $e_i^2$ in the law of friendship graphicable algebras $A(F_n)$ have a term more than in the case of star graphicable algebras $A(S_{2n})$. This term is
    $e_{i+1}$, if $i$ is even or $e_{i-1}$, if $i>1$ is odd, with the exception of $i = 2n+1$. \hfill $\Box$
    \end{lemma}

\vspace{0.15cm}

Similarly, a wheel graph can be obtained starting from the corresponding star graph by linking the vertex $e_i$ with $e_{i+1}$, for any subindex $i$ different from the corresponding with the central vertex. Thus, the following result is held

\begin{lemma}
 A wheel graphicable algebra $A(W_n)$ can be obtained from the corresponding star graphicable algebra $A(S_n)$ by adding the generators $e_{i-1}$ and $e_{i+1}$ in the expressions of each $e_i^2$ in the law of the star graphicable algebra, with $1\leq i \leq n-1$. \hfill $\Box$
\end{lemma}

\vspace{0.15cm}

As a conclusion, the following result is obtained

\begin{thm} \label{Main}
A star graphicable algebra $A(S_{2n})$ is a subalgebra of a friendship graphicable algebra $A(F_n)$ and this last one is, in turn, a subalgebra of a wheel graphicable algebra $A(W_{2n})$, con $n\leq 2$. \hfill $\Box$
\end{thm}

Recall that a subalgebra of a graphicable algebra is not graphicable in general. The previous study shows the existence of families of graphicable algebras in which this property is hereditary.

\subsection{Snark graphicable algebras}

\begin{figure}[h!]
\centering
\includegraphics[width=7.5cm]{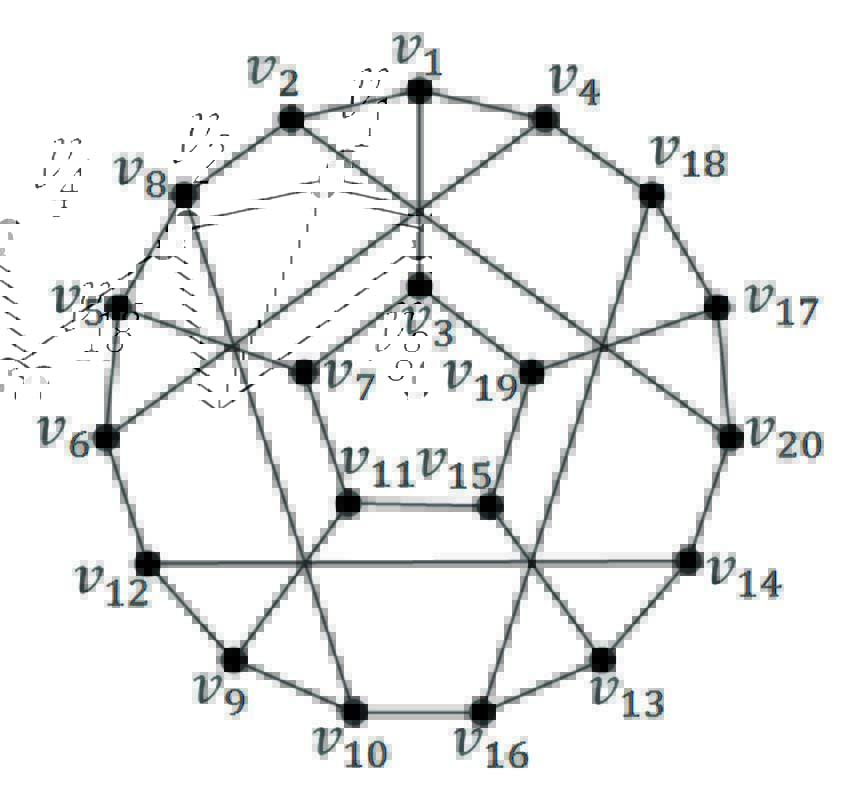}
\caption{Flower graph $J_5$.}
\label{Figura6}
\end{figure}

A {\em snark} is a connected, bridgeless cubic graph with chromatic index equal to 4. This last means that it is a graph in which every vertex has three neighbors, and the edges cannot be colored by only three colors without two edges of the same color meeting at a point. To avoid trivial cases, snark graphs are often restricted to have girth at least 5. For instance, a graph of this type is the so-called {\em flower} $J_5$, which has $20$ vertices. All its cycles have length greater than or equal to $5$. It is  represented in Figure \ref{Figura6}

We can associate this graph with a graphicable algebra in the following way

\begin{defi}
The $J_5$ graphicable algebra is that of dimension 20, with a set of generators $\{e_1, \ldots,e_{20}\}$ and law
$$
\begin{array}{ll}
e_i^2= e_{i+1}+e_{i+2}+e_{i+3}, & {\small {\it if \,\,\, i=1,5,9,13,17 \,\, (mod. \, 20)}}, \\
e_i^2= e_{i-1}+e_{i+6}+e_{i+18}, & {\small {\it if \,\,\,  i=2,6,10,14,18 \,\, (mod. \, 20)}}, \\
e_i^2= e_{i-3}+e_{i+2}+e_{i+14}, & {\small {\it if \,\,\,  i=4,8,12,16,20 \,\, (mod. \, 20)}}, \\
e_i^2= e_{i-2}+e_{i+4}+e_{i+16}, & {\small {\it if \,\,\,  i=3,7,11,15,19 \,\, (mod. \, 20)}}.
\end{array}
$$
\end{defi}

\vspace{0.15cm}

\begin{figure}[h!]
\centering
\includegraphics[width=7cm]{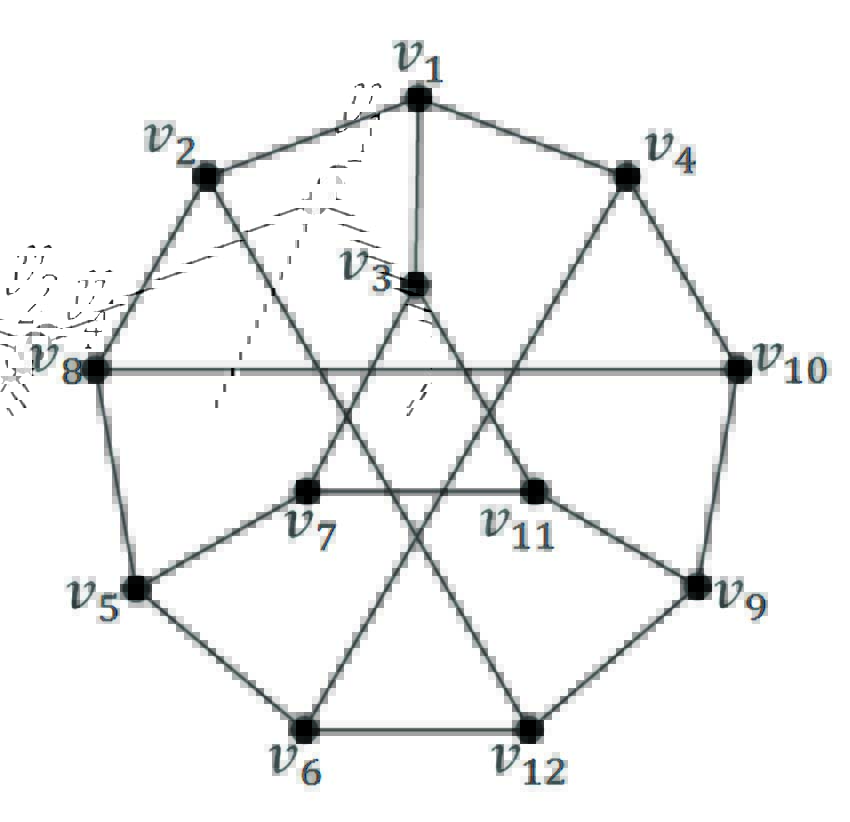}
\caption{Tietze graph.}
\label{Figura7}
\end{figure}

Another example of a graph of this family is the {\em Tietze graph}, which has $12$ vertices and $18$ edges. It is named after the mathematician
Heinrich Franz Friedrich Tietze (1880 - 1964). Like the Petersen graph, Tietze graph is maximally non-hamiltonian: it has no Hamiltonian cycle, but any two non-adjacent vertices can be connected by a Hamiltonian path (see Figure \ref{Figura7}).

\begin{defi}

The graphicable algebra of dimension 12, with a set of generators $\{e_1,\ldots,e_{12}\}$ and law
$$
\begin{array}{ll}
e_i^2= e_{i+1}+e_{i+2}+e_{i+3}, & {\small {\it if \,\,\,  i=1,5,9 \,\, (mod. \,  12)}}, \\
e_i^2= e_{i-1}+e_{i+6}+e_{i+10}, & {\small {\it if \,\,\,  i=2,6,10 \,\, (mod. \,  12)}}, \\
e_i^2= e_{i-3}+e_{i+2}+e_{i+6}, & {\small {\it if \,\,\, i=4,8,12 \,\, (mod. \,  12)}}, \\
e_i^2= e_{i-2}+e_{i+4}+e_{i+8}, & {\small {\it if \,\,\, i=3,7,11 \,\, (mod. \,  12)}}.
\end{array}
$$
\noindent is called the {\em Tietze graphicable algebra}.
\end{defi}
%\vspace{0.15cm}

Other snark graphs are the {\em Szekeres graph} (it was the fifth known snark, discovered by George Szekeres in 1973) which has 50 vertices and 75 edges and the {\em Watkins graph}, with the same numbers of vertices and edges than the former (Mark Watkins, 1964). In the same way, they can also be associated with the called {\em Szekeres graphicable algebra} and {\em Watkins graphicable algebra}, respectively. However, these algebras are not included here due to reasons of length.

\subsection{Graphicable algebras associated with generalized Petersen graphs}
\vspace{0.2cm}

The {\em generalized Petersen graphs} are a family of cubic graphs formed by connecting the vertices of a regular polygon to the corresponding vertices of a star polygon. This family, which includes the Petersen graph and generalizes one of the ways of constructing the Petersen graph, was introduced in 1950 by H.S.M. Coxeter and these graphs were given their name in 1969 by Mark Watkins.

Associated with these graphs, we can consider the following families of graphicable algebras in the following items.

\subsubsection{The Dürer graphicable algebra}

\begin{figure}[h!]
\centering
\includegraphics[width=7cm]{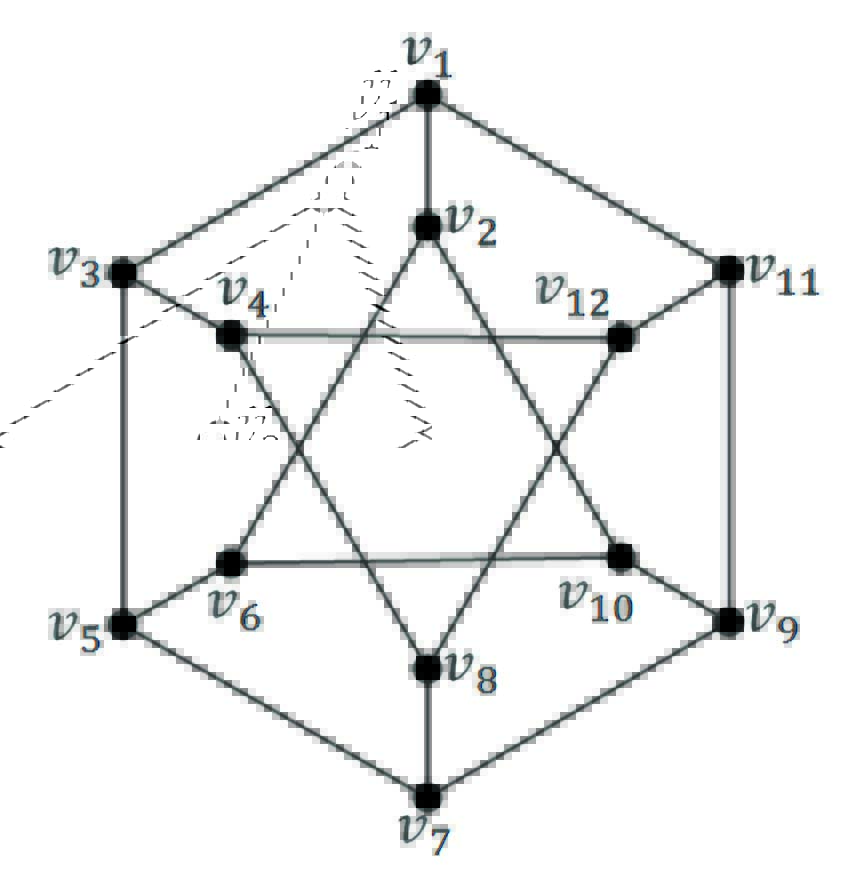}
\caption{Dürer graph.}
\label{Figura8}
\end{figure}

The {\em Dürer graph} is an undirected graph with 12 vertices and 18 edges, formed by the connection vertex to vertex between a regular hexagon and a six-point star. It is named after Albrecht Dürer (1471 - 1528), whose 1514 engraving {\it Melencolia I} includes a depiction of Dürer solid, a convex polyhedron having the Dürer graph as its skeleton. It is shown in Figure \ref{Figura8}.

The graphicable algebra associated with this graph is defined as follows

\begin{defi}

The {\em Dürer graphicable algebra} is that of dimension $12$, with the set of generators $\{e_1, \ldots,e_{12}\}$ and the law

$$
\begin{array}{ll}
  e_i^2=e_{i-2}+e_{i+1}+e_{i+2}, & {\small {\it if \,\,\, i \,\,es \,\, odd \, (mod. \,  12)}}, \\
  e_i^2=e_{i-1}+e_{i+4}+e_{i+8}, & {\small {\it if \,\,\, i \,\,es \,\, even \, (mod. \, 12)}}.
\end{array}
$$
\end{defi}

\subsubsection{The Möbius–Kantor graphicable algebra}

\begin{figure}[h!]
\centering
\includegraphics[width=7.5cm]{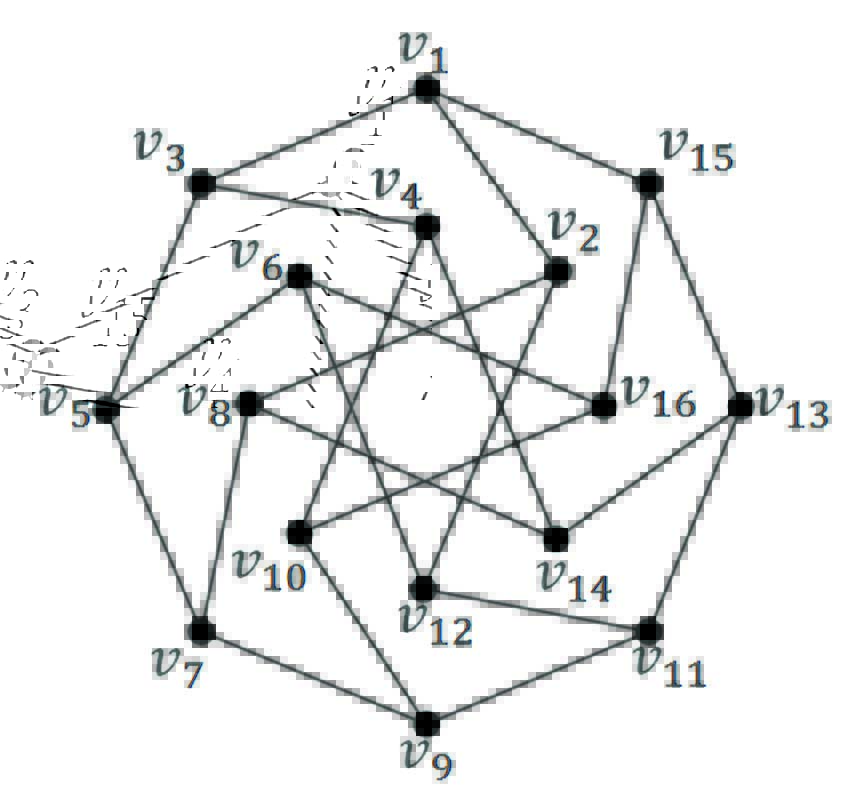}
\caption{Möbius-Kantor graph.}
\label{Figura9}
\end{figure}

The {\em Möbius–Kantor graph} is a non-planar, symmetric bipartite cubic graph with 16 vertices and 24 edges named after August Ferdinand Möbius (1790 - 1868) and Seligmann Kantor (German mathematician born in 1857 in Teplitz). It can be defined as the generalized Petersen graph of Figure \ref{Figura9}, that is, the one formed by the vertices of an octagon, connected to the vertices of an eight-point star such that it is shown in that figure.

Its corresponding associated graphicable algebra is defined as follows

\begin{defi}
The {\em Möbius-Kantor graphicable algebra} is that of dimension 16, with  a set of generators $\{e_1, \ldots,e_{16}\}$ and law

$$
\begin{array}{ll}
  e_i^2=e_{i-2}+e_{i+1}+e_{i+2}, & {\small {\it if \,\,\, i\,\, es\,\, odd \, (mod. \, 16)}}, \\
  e_i^2=e_{i-1}+e_{i+6}+e_{i+10}, & {\small {\it if \,\,\, i\,\, es\,\, even \, (mod. \, 16)}}.
\end{array}
$$
\end{defi}

\subsubsection{The Desargues graphicable algebra}

\begin{figure}[h!]
\centering
\includegraphics[width=7.5cm]{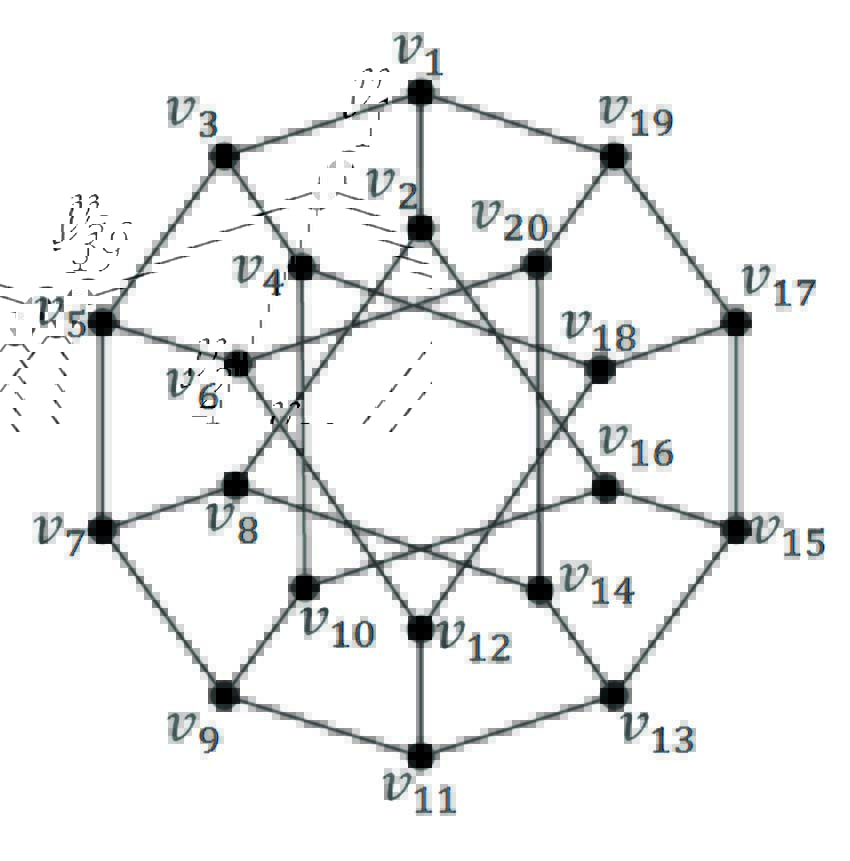}
\caption{Desargues graph.}
\label{Figura10}
\end{figure}

The {\em Desargues graph}, shown in Figure \ref{Figura10}, is a generalized Petersen graph formed by the vertices of an decagon, connected to the vertices of a ten-point star. It has 20 vertices and 30 edges and it is named after Gérard Desargues (1591 - 1661). It is Hamiltonian, 3-vertex-connected and 3-edge-connected.

%\vspace{0.15cm}

Its corresponding associated graphicable algebra is defined as follows

\begin{defi}
The {\em Desargues graphicable algebra} is that of dimension 20, with  a set of generators
$\{e_1,\ldots,e_{20}\}$ and law

$$
\begin{array}{ll}
  e_i^2=e_{i-2}+e_{i+1}+e_{i+2}, & {\small {\it if \, \, \, i\,\, es\,\, odd \, \, \, (mod. \,\, 20)}}, \\
  e_i^2=e_{i-1}+e_{i+6}+e_{i+14}, & {\small {\it if \, \, \, i\,\, es\,\, even \,\,\, (mod. \, \, 20).}}
\end{array}
$$
\end{defi}

\subsubsection{The Nauru graphicable algebra}

\begin{figure}[h!]
\centering
\includegraphics[width=7cm]{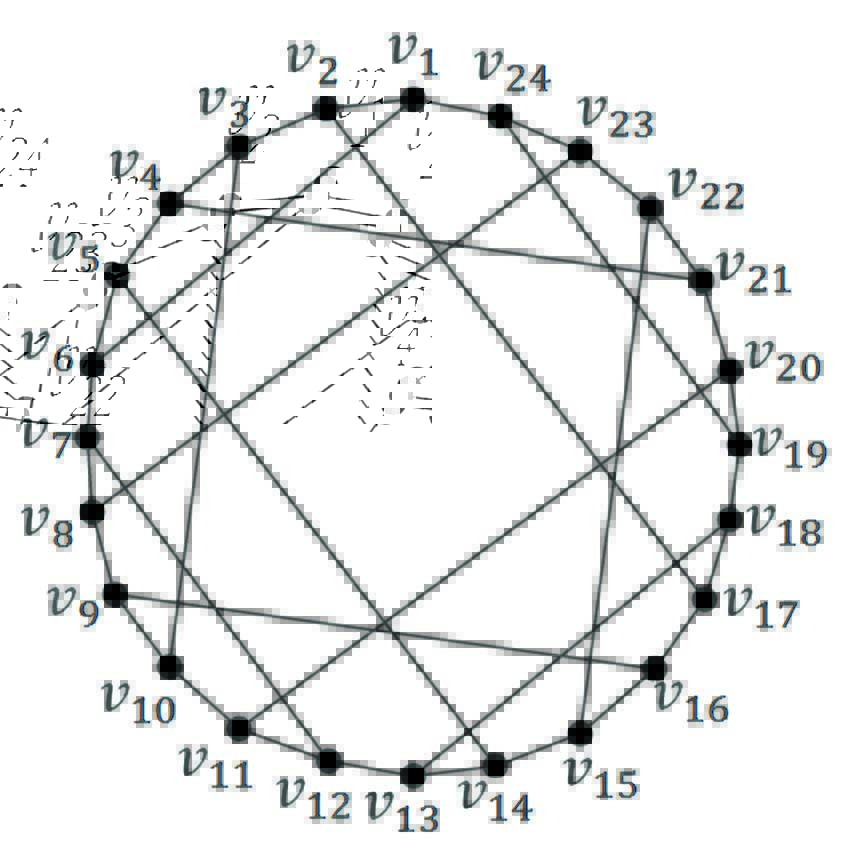}
\caption{Nauru graph.}
\label{Figura11}
\end{figure}

The {\em Nauru graph} is a bipartite cubic graph with 24 vertices and 36 edges. It is a generalized Petersen graph formed by the vertices of a dodecagon, connected to the vertices of an twelve-point star, as is shown in Figure
 \ref{Figura11}. It is 3-vertex-connected and 3-edge-connected.

%\vspace{0.4cm}

This graph was named by David Eppstein after the twelve-pointed star in the flag of Nauru, similar to the one that appears in the construction of the graph as a generalized Petersen graph. The Republic of Nauru, formerly known as Pleasant Island, is an island country in the Southern Hemisphere, whose length, 21 square kilometers makes it the smallest island nation in the world.

In the same way as before, we can associate this graph with a graphicable algebra, according to the following

\begin{defi}
The {\em Nauru graphicable algebra} is that of dimension 24, with  a set of generators
$\{e_1,\ldots,e_{24}\}$ and law

$$
\begin{array}{ll}
  e_i^2=e_{i-2}+e_{i+1}+e_{i+2}, & {\small {\it if \, \, \, i\,\, es\,\, odd \, \, \, (mod. \,\, 24)}}, \\
  e_i^2=e_{i-1}+e_{i+10}+e_{i+14}, & {\small {\it if \, \, \, i\,\, es\,\, even \,\,\, (mod. \, \, 24).}}
\end{array}
$$
\end{defi}

%$$
%\begin{array}{ll}
%  e_i^2=e_{i-1}+e_{i+1}+e_{i+5}, & {\small {\it if \, \, \, i=1,7,13,19 \,\, (mod. \,\, 24)}}, \\
%  e_i^2=e_{i-1}+e_{i+1}+e_{i+15}, & {\small {\it if \, \, \, i=2,8,14,20 \,\, (mod. \,\, 24)}}, \\
%  e_i^2=e_{i-1}+e_{i+1}+e_{i+7}, & {\small {\it if \, \, \, i=3,9,15,21 \,\, (mod. \,\, 24)}}, \\
%  e_i^2=e_{i-1}+e_{i+1}+e_{i+17}, & {\small {\it if \, \, \, i=4,10,16,22 \,\, (mod. \,\, 24)}}, \\
%  e_i^2=e_{i-1}+e_{i+1}+e_{i+9}, & {\small {\it if \, \, \, i=5,11,17,23 \,\, (mod. \,\, 24)}}, \\
%  e_i^2=e_{i-1}+e_{i+1}+e_{i+19}, & {\small {\it if \, \, \, i=6,12,18,24 \,\, (mod. \,\, 24).}}
%\end{array}
%$$

\subsubsection{Further results}

Apart from all the families of algebras here presented, authors are dealing with other new ones, such Platonic or fullerene graphicable algebras, associated with the families of Platonic and fullerene graphs, respectively.  Since all these graphs are planar, it would make sense to introduce the notion of {\em planar graphicable algebras} in general. They will be introduced in future work.

\end{document}